\documentclass{amsart}
\usepackage{amssymb}

\usepackage{graphicx,epstopdf}

\newtheorem{theorem}{Theorem}
\newtheorem{lemma}{Lemma}
\newtheorem{prop}{Proposition}

\begin{document}
\title
{ Periodic continued fractions and hyperelliptic curves}%
%
\author{M-P. Grosset}
\address{Department of Mathematical Sciences, Loughborough University, Loughborough, Leicestershire, LE11 3TU, UK}
         \email{M.Grosset@lboro.ac.uk}

         \author{A.P. Veselov}
\address{Department of Mathematical Sciences, Loughborough University, Loughborough, Leicestershire, LE11 3TU, UK and Landau Institute for Theoretical Physics, Moscow, Russia}

\email{A.P.Veselov@lboro.ac.uk}

\maketitle
\begin{abstract}

 We investigate when an algebraic function of the form $\phi(\lambda)=\frac{-B(\lambda)+\sqrt{R(\lambda)}}{A(\lambda)},$
 where $R(\lambda)$ is a  polynomial of odd degree $N=2g+1$ with coefficients in $\mathbb{C},$
 can be written as a periodic $\alpha$-fraction of the form
 $$\phi(\lambda)= [b_0;\overline{b_1, b_2, \dots, b_{N}}]_{\alpha} = b_0 + \frac{\lambda - \alpha_1}{b_1 +\frac{\lambda- \alpha_2}{b_2+_{\ddots  b_{N-1}+\frac{\lambda- \alpha_N}{b_N+  \frac{\lambda - \alpha_1}{b_1 +\frac{\lambda- \alpha_2}{b_2+_{\ddots  }}}}}}},$$
for some fixed sequence $\alpha_i.$  We show that this
problem has a natural answer given by the classical theory of hyperelliptic curves
and their Jacobi varieties.
We also consider pure periodic $\alpha$-fraction expansions corresponding to the special case when $b_N=b_0.$ 
\end{abstract}

\section{Introduction}
Consider the following continued fraction, which we will call {\it
$\alpha$-fractions}:
\begin{equation}
\label{Cont_frac1} \phi = b_0 + \frac{\lambda - \alpha_1}{b_1
+\frac{\lambda- \alpha_2}{b_2+_{\ddots}}} = [b_0,b_1,\dots,
]_{\alpha},
\end{equation}
where $\alpha = (\alpha_i), \alpha_i \in \mathbb{C}$ is a given
sequence, $b_i$ are arbitrary complex numbers, $\lambda$ is a
formal parameter. In this paper we will consider a special case of
{\it $N$-periodic $\alpha$-fractions}, when the
sequences $\alpha_i$ and $b_i$ are periodic with period $N:$ $$
\alpha_{i+N}=\alpha_i,\,
 b_{i+N}=b_{i}
$$ for all $i\geq 1:$
\begin{equation}
\label{per}\phi=[b_0;\overline{b_1, b_2, \dots, b_{N}}]_{\alpha}.
\end{equation}
In the particular case when $b_N=b_0$ we have
$\phi=[\overline{b_0, b_1, \dots, b_{N-1}}]_{\alpha},$ which will
be called  a {\it pure N-periodic $\alpha$-fraction.}

This kind of fractions naturally appear in the theory of
integrable systems, in particular in the theory of periodic
dressing chain \cite{VS}, but to the best of our knowledge has not
been studied so far. We were partly inspired by our recent
discussions with Vassilis Papageorgiou on the discrete KdV
equation where such continued fractions appear as well
\cite{GPV}.

Because of periodicity we can write formally (\ref{per}) as
$$\phi= b_0 + \frac{\lambda - \alpha_1}{b_1 +\frac{\lambda-
\alpha_2}{b_2+_{\ddots b_{N-1}+\frac{\lambda-
\alpha_N}{b_N-b_0+\phi}}}},$$ which implies a quadratic relation
\begin{equation}
\label{Quad_eq} A(\lambda) \phi^2 + 2 B(\lambda) \phi + C(\lambda)
=0,
\end{equation}
where $A, B, C$ are certain polynomials in $\lambda$ with
coefficients polynomially depending on $b_i.$ Thus to any periodic
$\alpha$-fraction (\ref{per}) corresponds an algebraic
function
\begin{equation}
\label{Algebraic_fct}
\phi(\lambda)=\frac{-B(\lambda)+\sqrt{R(\lambda)}}{A(\lambda)},
\end{equation}
 where
\begin{equation}
\label{R} R(\lambda)=B(\lambda)^2 - A(\lambda)C(\lambda)
\end{equation}
is the discriminant of (\ref{Quad_eq}). In that case we will say
that (\ref{per}) is a {\it periodic $\alpha$-fraction
expansion} of the algebraic function (\ref{Algebraic_fct}) from the
hyperelliptic extension $\mathbb{C}(\lambda, \sqrt{R(\lambda)})$
of the field of rational functions $\mathbb{C}(\lambda)$. We leave the question of convergence
aside concentrating on algebraic and geometric aspects of the problem.

We will discuss the following three main questions.

{\bf Question 1.} Which algebraic functions (\ref{Algebraic_fct})
admit $N$-periodic $\alpha$-fraction expansions ?

{\bf Question 2.} How many such expansions may exist for a given
algebraic function (\ref{Algebraic_fct}) and how to find them ?

{\bf Question 3.} What is the geometry of the set of functions
(\ref{Algebraic_fct}) from given hyperelliptic extension (i.e.
with fixed $R$), which admit periodic $\alpha$-fraction
expansions?

The answers depend on the parity of $N.$ In this paper we will
restrict our study to the case of odd period $N= 2g+1,$ which is
the most interesting one (cf. \cite{VS}). We will also assume that
all the parameters $\alpha_i$ are distinct.

Note that when $N$ (which is also the degree of
$R(\lambda)$) is even, one can consider the usual continued
fraction expansions going back to Abel and Chebyshev who
discovered their relation with the classical problem of integration
in elementary functions (see a nicely written paper by van der
Poorten and Tran \cite{vanPoorten} for details). These expansions
are more natural, but can not be used in the odd degree case.

We would like to mention also that in the classical number-theoretic version
the answer to Question 1 is due to Galois, who proved that a quadratic irrationality
$\xi = p + q\sqrt{d}$,  $p$ and $q$ rational numbers, $d$ integer, has a pure periodic continued fraction expansion
$\xi = [\overline{a_0,\dots, a_k}]$ if and only if $\xi$ is larger than 1 and its conjugate $\bar \xi = p - q\sqrt{d}$ lies between $-1$ and $0$ (see e.g. \cite{LeVeque}). Note also that the periodic continued
fractions of the form $[a_0; \overline{a_1,\dots, a_k}]$ also appear naturally in number theory as expansions of $\sqrt{d}$ (see \cite{LeVeque}). 

To explain our main results let us introduce the polynomial
\begin{equation}
\label{Alpha} \mathfrak{A}(\lambda)= \prod_{i=1}^{N} (\lambda-
\alpha_i)
\end{equation}
and call a polynomial $R$ of degree $N$ {\it $\alpha$-admissible}
if
\begin{equation}
\label{R} R(\lambda) = S^2(\lambda) + \mathfrak{A} (\lambda)
\end{equation}
for some polynomial $S(\lambda)$ of degree $g$ or less,
where as before $N = 2g+1.$
We will call a polynomial {\it monic}
 if its highest coefficient is $1$ and {\it anti-monic} if it is equal to $-1.$
Note that $\alpha$-admissible polynomials $R$ are automatically
monic.

\begin{theorem} \label{theorem1}
The algebraic functions $\phi(\lambda)$ admitting an $N$-periodic
$\alpha$-fraction expansion have the form
(\ref{Quad_eq}, \ref{Algebraic_fct}) with the polynomials $A,\, B,\,
C$ satisfying the following conditions:
\begin{enumerate}
\item $deg \, B \leq g, \quad  A(\lambda)$ and $C(\lambda)$ are monic
and anti-monic polynomials of degree $g$ and $g+1$ respectively
\item the discriminant $R(\lambda)= B^2 - AC$ is $\alpha$-admissible.
\end{enumerate}
Conversely, for an open dense subset of such triples $(A, B, C)$
the corresponding function (\ref{Algebraic_fct}) has exactly two
$N$-periodic $\alpha$-fraction expansions. The
corresponding $b_i$ are rational functions of both coefficients of
$A,\, B,\,C$ and parameters $\alpha_i$ and can be found by an
effective matrix factorisation procedure.

In the pure periodic case the only additional requirement is
\begin{equation}
\label{add} C(\alpha_N) = 0,
\end{equation}
 under which the pure
periodic $\alpha$-fraction expansion is generically unique.
\end{theorem}

We will call $(A, B, C)$ satisfying the conditions (1), (2) of
Theorem 1 the {\it $\alpha$-triples}. Note that these conditions are
invariant under any permutation of the parameters $\alpha_i.$
In fact there is a natural birational action of the direct product
$G = \mathbb{Z}_2 \times S_N$ on the set of periodic continued
$\alpha$-fraction expansions, where the generator $\varepsilon$ of
$\mathbb{Z}_2$ is acting simply by swapping two different
$\alpha$-fraction expansions given by Theorem 1. Our next result
describes this action explicitly.

Let us introduce the following permutation $\pi \in S_N$, which
reverses the order $\alpha_1, \alpha_2, \dots, \alpha_{N-1},\, \alpha_N$ to $\alpha_N,\, \alpha_{N-1}, \dots, \alpha_2,
\alpha_1$ and the involutions $\sigma_k$ swapping $\alpha_k$ and $\alpha_{k+1}$, where $k
= 1, \dots, N-1$.

We will show that $\sigma_k$ is acting on $b = (b_i), \, i=0,
\dots, N$ with $b_k \neq 0$ as follows:
\begin{equation}
\label{actionsigma} \tilde b_{k-1} = b_{k-1}+
\frac{\alpha_{k+1}-\alpha_k}{b_k} ,\quad \tilde b_{k+1}= b_{k+1}
-\frac{\alpha_{k+1}-\alpha_k}{b_k},
\end{equation}
the rest of $b_i$ remain the same. This determines the action of
the symmetric group $S_N$ since $\sigma_k$ generate it. To
describe the action of $\mathbb{Z}_2$ it is enough to describe the
action of the involution $\varepsilon \pi \in G,$ which turns out to be
quite simple:
\begin{equation}
\label{actionipi} \tilde b_{j} = - b_{N-j}, \, j = 1, \dots, N-1
\quad \tilde b_{0} = b_0 - b_{N}, \quad \tilde b_{N} = - b_{N}.
\end{equation}

\begin{theorem} \label{theorem2}
The formulae (\ref{actionsigma}) and (\ref{actionipi}) define a
birational action of the group $G = \mathbb{Z}_2 \times S_N$ on
the set of $N$-periodic $\alpha$-fractions. Its orbits
consist of all $2 N!$ possible periodic $\alpha$-fraction expansions for a given $\alpha$-triple $(A, B,
C)$ and any permutation of the parameters $\alpha_i.$

In the pure periodic case the symmetry group is broken down to
$S_{N-1}$ generated by $\sigma_k$ with $k=1, \dots, N-2$ given by
(\ref{actionsigma}).
\end{theorem}

Let us fix now the $\alpha$-admissible polynomial $R(\lambda)$ with
distinct roots. We would like to describe the geometry of the set of elements from
the hyperelliptic extension field $\mathbb{C}(\lambda,
\sqrt{R(\lambda)})$ which have a pure periodic
$\alpha$-fraction expansion. For the standard material from the algebraic geometry of the curves we refer to the classical Griffiths-Harris book \cite{GH}.

Consider the hyperelliptic curve
$\Gamma_R$ given by the equation
\begin{equation}
\label{Gamma} \mu^2 = R(\lambda).
\end{equation}
The curve $\Gamma_R$ consists of the affine part $\Gamma_R^{aff},$ corresponding to the "finite" solutions of (\ref{Gamma}), and the "infinity" point, which we will denote as $P_{\infty}.$
Since the roots of $R$ are distinct it is non-singular and has genus $g.$
Consider $g$ points $P_1, \dots, P_g$ of $\Gamma_R^{aff}$ and call the corresponding divisor
$D = P_1 + \dots + P_g$ {\it non-special} if 
\begin{equation}
\label{nonspec}
P_i \neq \tau(P_j)
\end{equation}
for any $i \neq j$, where $\tau$ is
the {\it hyperelliptic involution}: $$\tau (\mu,\lambda) =
(-\mu,\lambda).$$
Non-special divisors have the property that the linear space $L(D)$ of all meromorphic functions on $\Gamma_R$ having the poles at $D$ of order less than or equal to 1 has dimension 1, which means that it consists only of constant functions. The corresponding linear space $L(D + P_{\infty})$ has dimension 2, so there exists a  non-constant function $f \in L(D + P_{\infty})$ with additional pole at infinity. These functions are in a way "least singular" among the "generic" meromorphic functions on $\Gamma_R$ in the sense that any such function can not have less than $g+1$ poles (see \cite{GH}).

Now define the {\it affine Jacobi variety} $J(\Gamma_R)^{aff}$
as the set of  positive non-special divisors $D = P_1 + \dots + P_g,\,
P_1, \dots, P_g \in \Gamma_R^{aff}.$

Let $M_R^{\alpha}$ be an affine variety of $\alpha$-triples of
polynomials $(A,\, B,\, C)$ with given discriminant $R(\lambda),$
and $P_R^{\alpha}$ be its subvariety given by the
additional condition (\ref{add}).

\begin{theorem} \label{theorem3}
There exists a bijection between the set $M_R^{\alpha}$ and the
extended affine Jacobi variety $J(\Gamma_R)^{aff} \times
\mathbb{C}.$ The corresponding algebraic functions (\ref{Algebraic_fct})
 can be characterised as meromorphic functions $\phi \in L(D + P_{\infty})$ on $\Gamma_R$ with non-special pole divisor
$D+P_{\infty}$ and asymptotic $\phi \sim \sqrt{\lambda}$ at infinity.

In the pure periodic case under the assumption that $R(\alpha_N) \neq
0$ there exists a natural $2:1$ map from the set $P_R^{\alpha}$ to $J(\Gamma_R)^{aff}.$ The corresponding $\phi$ from $L(D + P_{\infty})$
are fixed by the condition that one of two values of $\phi(\alpha_N)$ is zero.
\end{theorem}

The proof is based on the classical description of the Jacobi
variety due to Jacobi himself \cite{Jacobi} (see also Mumford
\cite{Mum}). 

We see that the affine space $\mathbb{C}^N$ of all functions
(\ref{Algebraic_fct}) having periodic $\alpha$-fraction
expansion is birationally equivalent to the double covering of the
bundle of the extended affine Jacobians of $\alpha$-admissible
hyperelliptic curves. This is a version of the well-known result
by Dubrovin and Novikov \cite{DN} who were the first to apply the
theory of the KdV equation to the problems of algebraic geometry.

\section{Periodic $\alpha$-fractions \label{section periodic  continued fractions}}

Consider the  $N$-periodic $\alpha$-fraction of period
$N=2g+1$ \\ $$\phi=  b_0 + \frac{\lambda - \alpha_1}{b_1
+\frac{\lambda- \alpha_2}{b_2+_{\ddots  b_{N-1}+\frac{\lambda-
\alpha_N}{b_N + \frac{\lambda - \alpha_1}{b_1 +\frac{\lambda-
\alpha_2}{b_2+_{\ddots }}}}}}} =[b_0; \overline {b_1, \dots,
b_{N-1},b_N}]_{\alpha}$$

We see that (at least formally) $\phi$ is the fixed point of the
fractional linear transformation
\begin{equation}
\label{fraclin} s(\phi)=b_0 + \frac{a_1}{b_1
+\frac{a_2}{b_2+_{\ddots b_{N-1}+\frac{a_N}{b_N^*+\phi}}}}
\end{equation}
with $b_N^*=b_N-b_0$ and $a_k=\lambda-\alpha_k, \, k=1,..,N.$  The
function $s(\phi)$ can be written as $s(\phi)=\frac{P_{N-1}\phi +
P_{N}}{Q_{N-1}\phi + Q_{N}},$ where the quantities $P_{k},\,Q_{k}$
are determined by the standard recurrence relations (see e.g.
\cite{Wall}, page 14): $$P_{-1} =1,\quad Q_{-1}=0,$$$$
P_0=b_0,\quad Q_0=1,$$
\begin{equation}
\label{recur}
P_{k+1}=b_{k+1}P_k+a_{k+1}P_{k-1},~Q_{k+1}=b_{k+1}Q_k+a_{k+1}Q_{k-1}
\end{equation}
for $k < N-1$ and $$P_{N}=b_{N}^*P_{N-1}+a_{N}P_{N-2}, \quad
Q_{N}=b_{N}^*Q_{N-1}+a_{N}Q_{N-2}. $$

Thus we have
\begin{equation}
\label{phi} \phi=\frac{P_{N-1}\phi + P_{N}}{Q_{N-1}\phi +Q_{N}},
\end{equation}
which can be written as a quadratic equation
\begin{equation}
\label{Quadeq2} Q_{N-1} \phi^2 +( Q_{N} -P_{N-1})\phi -P_{N} =0.
\end{equation}
It is easy to see from the recurrence relations that $P_k$ and
$Q_k$ are polynomials in $\lambda$ of the form $$
P_{2k}=(b_0+b_2+...+b_{2k})\lambda^k+ \dots, \quad
P_{2k+1}=\lambda^{k+1}+  \dots,$$ $$ Q_{2k}=\lambda^k+
 \dots, \quad Q_{2k-1}=(b_1+b_3+...+b_{2k-1})\lambda^{k-1}+  \dots,$$
for $k \leq g$ and
$Q_{2g+1}=(b_1+b_3+...+b_{2g-1}+b_{2g+1}-b_0)\lambda^g+ \dots,$ where the dots denote the lower degree terms. Since from
(\ref{Quadeq2})
\begin{equation}
\label{pq} A(\lambda)= Q_{N-1} (\lambda),\quad B(\lambda)
=\frac{1}{2}( Q_{N}(\lambda) -P_{N-1}(\lambda)), \quad C(\lambda)=
- P_{N}(\lambda),
\end{equation}
the polynomial $A$ is monic of degree $g$, $C$ is anti-monic of
degree $g+1$ and $B$ has degree $g$ or less with the highest
term $\beta\lambda^g,$ where $$\beta =  -b_0
+\frac{1}{2}\sum_{k=1}^{N} (-1)^{k+1}b_k.$$

Let us show now that the discriminant $R = B^2 -AC$ is
$\alpha$-admissible. We have $$R=\frac{1}{4}( Q_{N} -P_{N-1})^2
+P_{N}Q_{N-1} = \frac{1}{4}(P_{N-1}+ Q_{N} )^2
+P_{N}Q_{N-1}-P_{N-1} Q_{N}.$$ We claim that $$
P_{N}Q_{N-1}-P_{N-1} Q_{N} = \prod_{i=1}^{N} (\lambda-
\alpha_i).$$ Indeed, the determinant $$\left | \begin {array}{cc}
P_{N} & P_{N-1}
\\Q_{N}&Q_{N-1} \end{array} \right |= \left | \begin {array}{cc}
b_N P_{N-1} + a_NP_{N-2} &P_{N-1} \\ b_NQ_{N-1}+
a_{N}Q_{N-2}&Q_{N-1}  \end{array} \right |=$$ $$- a_N \left |
\begin {array}{cc} P_{N-1} & P_{N-2}  \\Q_{N-1}&Q_{N-2}
\end{array} \right |=\dots= (-1)^{N} a_N a_{N-1} \dots  a_1 \left
| \begin {array}{cc} b_0 &1\\ 1 &0\end{array} \right |.$$ Since
$N$ is odd, $$\left | \begin {array}{cc} P_{N} & P_{N-1}
\\Q_{N}&Q_{N-1} \end{array} \right |=a_N a_{N-1} \dots  a_1 =
\prod_{i=1}^{N} (\lambda- \alpha_i)= \mathfrak{A}.$$

Now by taking $S(\lambda) =\frac{P_{N-1}+ Q_{N} }{2},$ which is
a polynomial of degree $g$ or less, we see that
$R(\lambda) = S^2 + \mathfrak{A},$ so $R$ is $\alpha$-admissible.
This proves the first part of Theorem 1 in the periodic case.

To prove the second part let us introduce the following matrix
\begin{equation}
\label{mat1} M(\lambda)= \left [ \begin {array}{cc}1 &b_0 \\ 0&1
\end{array} \right ] \left [ \begin {array}{cc}0 & \lambda-
\alpha_1 \\1&b_1 \end{array} \right ] \dots \left [ \begin
{array}{cc} 0 & \lambda- \alpha_N \\1&b_N^*\end{array} \right ]  ,
\end{equation}
with $b_N^*=b_N-b_0.$ One can check that it can be rewritten also
as
\begin{equation}
\label{mat2} M = \left [ \begin {array}{cc} b_0 & \lambda-
\alpha_1
\\1&0
\end{array} \right ] \dots \left [ \begin {array}{cc} b_{N-1} &
\lambda- \alpha_N \\1&0 \end{array} \right ] \left [ \begin
{array}{cc}1 & b_{N}^*\\0&1 \end{array} \right ] .
\end{equation}

The following Lemma explains its importance for our problem.

\begin{lemma} \label{matrix1}
Vector $\left( \begin {array}{c} \phi \\ 1 \end{array} \right )$
with $\phi=[b_0; \overline {b_1, \dots, b_{N-1},b_N}]_{\alpha}$ is
an eigenvector of the matrix $M(\lambda).$
\end{lemma}

The proof follows from the fact that $\phi$ is the fixed point of
the fractional linear transformation (\ref{fraclin}). The product
of matrices (\ref{mat1}) corresponds to the representation of
$s(\phi)$ as a superposition $s_0\circ s_1 \circ \dots \circ
s_N(\phi),$ where $s_0(\phi) =b_0+\phi,
s_k(\phi)=\frac{\lambda-\alpha_k}{b_k+\phi}$ for $k=1,2,..., N-1$
and $ s_N(\phi)=\frac{\lambda-\alpha_N}{b_N^*+\phi}.$

Let $T(\lambda)= \frac{1}{2} tr \, M$ be half of the trace of the
matrix $M(\lambda),$ which is a polynomial of degree $g$ or less.
Note that the determinant of $M$ is equal to
$-\mathfrak{A} = -\prod_{i=1}^{N} (\lambda- \alpha_i)$ as it
follows immediately from (\ref{mat1}).

\begin{lemma} \label{matrix2}
The matrix (\ref{mat1}) has the form
\begin{equation}
\label{matform} M(\lambda)=\left [ \begin
{array}{cc}T(\lambda)-B(\lambda) & -C(\lambda)
\\A(\lambda)&T(\lambda)+B(\lambda) \end{array} \right ],
\end{equation}
where $(A, B, C)$ is the $\alpha$-triple of polynomials corresponding to
$\phi.$ The discriminant $R = B^2 - AC$ equals to $T^2 +
\mathfrak{A}.$
\end{lemma}

Indeed $$M(\lambda)=\left [ \begin {array}{cc} P_{N-1}(\lambda) &
P_N(\lambda)
\\Q_{N-1}(\lambda)& Q_N(\lambda) \end{array} \right ],$$
where $P_k, \, Q_k$ are defined above by (\ref{recur}). Now the
first claim follows from the relations (\ref{pq}). Taking the
determinant of both sides of (\ref{matform}) we have
$-\mathfrak{A} = T^2 - B^2 + AC$, which implies $B^2 -AC = T^2 +
\mathfrak{A}.$

Now we need the following result about the factorisation of such
matrices. This kind of problems often appears in the theory of
discrete integrable systems (see \cite{MV} and \cite{AV}).

\begin{prop} \label{factorisation}
Let $M(\lambda)$ be a polynomial matrix of the form (\ref{matform}),
where $A$ is a monic polynomial of degree $g$, $C$ is an
anti-monic polynomial of degree $g+1$, $T$ and $B$ are polynomials of degree
$g$ or less. Assume also that $\det M(\lambda) = -\prod_{i=1}^{N}
(\lambda- \alpha_i).$ Then for an open dense set of such $M$ there
exists a unique factorisation of the form $$M(\lambda)=\left [
\begin {array}{cc} b_0 & \lambda- \alpha_1 \\1&0 \end{array}
\right ] \dots \left [ \begin {array}{cc} b_{N-1} & \lambda-
\alpha_N \\1&0
\end{array} \right ] \left [ \begin {array}{cc}1 & b_{N}-b_0\\0&1
\end{array} \right ].$$
\end{prop}

The proof is actually effective. We describe the procedure which
allows to find $b_i$ uniquely assuming at the beginning that the
factorisation exists.

Consider the transpose $M^T$ of the matrix $M.$ For $\lambda=
\alpha_1$ the matrix $M^T(\lambda)$ is degenerate (since $\det
M^T(\lambda) = \det M(\lambda) = -\prod_{i=1}^{N} (\lambda-
\alpha_i)$). Find the {\it null-vector} $ e_1 = {\left (
\begin {array}{c} x_1
\\ y_1
\end{array} \right ) }$ of $ M^T(\alpha_1)$, which is by definition any non-zero vector such that
\begin{equation}
\label{e1} M^T(\alpha_1) e_1=0,
\end{equation}
or explicitly $$ \left [
\begin {array}{cc}T(\alpha_1)-B(\alpha_1) & C(\alpha_1)  \\
A(\alpha_1)&T(\alpha_1)+B(\alpha_1) \end{array} \right ] {\left (
\begin {array}{c} x_1
\\ y_1
\end{array} \right ) }= {\left (
\begin {array}{c} 0
\\ 0
\end{array} \right ). }
$$
 It must satisfy the relation $$\left [
\begin {array}{cc} b_0 & 1 \\0&0 \end{array}
\right ] {\left (
\begin {array}{c} x_1
\\ y_1
\end{array} \right ) } = {\left (
\begin {array}{c} 0
\\ 0
\end{array} \right ) }
$$ since all other factors are non-degenerate when $\lambda =
\alpha_1.$ This determines $b_0$ uniquely as
\begin{equation}
 \label{b0}
b_0= \frac{T(\alpha_1)-B(\alpha_1)}{A(\alpha_1)}= \frac{C(\alpha_1)}{T(\alpha_1)+B(\alpha_1)}.
\end{equation}

Consider now the matrix $ M_1 = {\left [ \begin {array}{cc} b_0 &
\lambda- \alpha
\\1&0
\end{array} \right ] }^{-1}M(\lambda).$
It is polynomial in $\lambda$ because of the following elementary

\begin{lemma}
Let $M$ be the polynomial matrix, $\lambda=\alpha$ be a simple
root of its determinant and  $e = {\left (
\begin {array}{c} 1  \\ -b
\end{array} \right ) }$ be a null vector of $M^T(\alpha).$ Then the
matrix $ {\left [ \begin {array}{cc} b & \lambda- \alpha \\1&0
\end{array} \right ] }^{-1}M(\lambda)$
is polynomial.
\end{lemma}

Indeed, let $M= \left [ \begin {array}{cc} X(\lambda) & Y(\lambda)
\\ Z(\lambda)&W(\lambda) \end{array} \right ]$ then  $$ {\left [
\begin {array}{cc} b & \lambda- \alpha
\\1&0
\end{array} \right ] }^{-1}M(\lambda)=  {\left [ \begin
{array}{cc} Z(\lambda) &W(\lambda) \\
\frac{X(\lambda)-bZ(\lambda)}{ \lambda-
\alpha}&\frac{Y(\lambda)-bW(\lambda)}{ \lambda- \alpha}
\end{array} \right ] }.$$
From $M^T(\alpha) {\left [ \begin {array}{c} 1  \\-b \end{array}
\right ] }=0$ it follows that $\lambda=\alpha$ is a root of the
polynomials $X(\lambda)-bZ(\lambda)$ and $Y(\lambda)-bW(\lambda).$
Therefore these polynomials are divisible by $\lambda-\alpha,$
which proves the claim.

Repeat now the procedure by taking $\lambda = \alpha_2$ and so on.
After $N$ steps we will come to a polynomial matrix
$$M_N(\lambda)= {\left [
\begin {array}{cc} b_{N-1} & \lambda- \alpha_{N} \\1&0 \end{array}
\right ] }^{-1}\times \dots  \times{\left [ \begin {array}{cc} b_0
& \lambda- \alpha_1 \\1&0 \end{array} \right ] }^{-1}M(\lambda)$$
with determinant 1. To complete the proof of Proposition
\ref{factorisation} we need to show that $M_N$ is of the form
$\left [ \begin {array}{cc}1 & b_{N} ^*\\0&1 \end{array} \right
].$

Recall that the matrix $M(\lambda)$ is of the form $ \left [
\begin {array}{lr}a_0 \lambda^g+\dots & \lambda^{g+1}+\dots \\
\lambda^{g}+\dots &d_0 \lambda^g+\dots,\end{array} \right ], $
where the dots mean terms of lower degree, and the coefficients
$a_0$ and $d_0$ may be zero. It is easy to show that the matrix
$M_2(\lambda)={\left [ \begin {array}{cc} b_1 & \lambda- \alpha_2
\\1&0 \end{array} \right ] }^{-1}{\left [ \begin {array}{cc} b_0 &
\lambda- \alpha_1 \\1&0 \end{array} \right ] }^{-1}M(\lambda)$ is
of the form $ \left [ \begin {array}{lr}a_2\lambda^{g-1}+\dots &
\lambda^g+\dots \\ \lambda^{g-1}+\dots&d_2 \lambda^{g-1}+\dots
\end{array} \right ]$
and by induction  $M_{2k}(\lambda)$ is of the form
 $ \left [ \begin {array}{lr}a_k\lambda^{ g-k}+\dots & \lambda^{g-k+1}+\dots \\ \lambda^{g-k}+\dots&d_k\lambda^{ g-k}+\dots \end{array} \right ]. $
Therefore the matrix $M_{N-1}$ is of the form $\left [ \begin
{array}{cc} a_g & \lambda+c \\1&d_g \end{array} \right ] $ where
$a_g,c,d_g$ are  constant. The matrix $M_{N}= {\left [
\begin {array}{cc} b_{N-1} & \lambda- \alpha_{N} \\1&0 \end{array}
\right ] }^{-1} M_{N-1}$ equals to $\left [ \begin {array}{cc} 1
&d_g \\ \frac{a_g-b_{N-1}}{ \lambda- \alpha_{N}}&\frac{\lambda +
c-b_{N-1}d_g}{ \lambda- \alpha_{N}} \end{array} \right ].$ Since
$M_N$ is a polynomial matrix, we have $a_g=b_{N-1}$ and $
b_{N-1}d_g -c =\alpha_N .$ Thus $M_N$ has the required form.

We see that the procedure will not work only if at some stage the
first component of the null vector of $M_k^T(\alpha_{k+1})$
vanishes. Clearly this happens only for a closed algebraic subset
of codimension 1, so for generic triples $(A,B,C)$ the matrix
decomposition exists and is unique. This completes the proof of
Proposition \ref{factorisation}.

Now we are ready to finish the proof of Theorem 1 in the periodic
case. Let $(A, B, C)$ be an $\alpha$-triple, then by definition there
exists a polynomial matrix $S$ of degree $g$ or less 
such that the discriminant $R = B^2 -AC$ is equal to $S^2 +
\mathfrak{A}.$ Clearly the polynomial $S$ is defined up to a sign.
Consider two corresponding matrices $M$ given by (\ref{matform})
with $T(\lambda) = \pm S(\lambda).$ Each of them generically has a
unique factorisation given by Proposition 1. One can easily check
that this gives two $N$-periodic $\alpha$-fraction
representations of the corresponding function
$\phi(\lambda)=\frac{-B(\lambda)+\sqrt{R(\lambda)}}{A(\lambda)}$
and thus completes the proof in this case.

\section{Pure periodic $\alpha$-fractions \label{section continued fractions}}

Let now $\phi=[\overline{b_0, b_1, \dots, b_{N-1}}]_{\alpha}$ be a
pure periodic $\alpha$-fraction. This is a particular
case of the previous situation with $b_0 = b_N.$ But since the
corresponding $b_N^* = b_0 - b_N = 0,$ this case is actually
degenerate and needs a special consideration.

First of all as before $\phi$ satisfies the relation $$
\phi=\frac{P_{N-1}\phi + P_{N}}{Q_{N-1}\phi +Q_{N}},$$ where $P_k,
\, Q_k$ satisfy the relations (\ref{recur}), but now because
$b_N^* = 0$ we have $$P_N = (\lambda - \alpha_N) P_{N-2}, \quad
Q_N = (\lambda - \alpha_N) Q_{N-2}.$$ Now from (\ref{pq}) we have
$$ A(\lambda)= Q_{N-1} (\lambda),\, B(\lambda) = \frac{1}{2}(
(\lambda - \alpha_N) Q_{N-2}(\lambda) - P_{N-1}(\lambda)), \,
C(\lambda)
 = - (\lambda -
\alpha_N) P_{N-2}. $$ Since $Q_{N-1}$ and $P_{N-2}$ are monic this
shows that $A, B, C$ satisfy the property (1) of Theorem 1 with the
additional condition $C(\alpha_N)=0$. The proof of the second
property ($\alpha$-admissibility of $R$) goes unchanged.

Now as in the previous case in order to find the pure periodic
$\alpha$-fraction one should factorise the matrix
$$M(\lambda)=\left [ \begin {array}{cc}T(\lambda)-B(\lambda) & -
C(\lambda)
\\A(\lambda)&T(\lambda)+B(\lambda) \end{array} \right ]$$ as the
product $$\left [ \begin {array}{cc} b_0 & \lambda- \alpha_1 \\1&0
\end{array} \right ] \dots \left [ \begin {array}{cc}
b_{N-1} & \lambda- \alpha_N \\1&0 \end{array} \right ].$$ The main
difference is that in the pure periodic case the trace
$T(\lambda)$ of the matrix $M(\lambda),$ which was known before
only up to a sign, now is determined uniquely by the condition
\begin{equation}
\label{tr} T(\alpha_N) = - B(\alpha_N).
\end{equation}
Indeed $T(\alpha_N) =\frac{1}{2}(P_{N-1}(\alpha_N)+ Q_N(\alpha_N))
= \frac{1}{2}(P_{N-1}(\alpha_N)) = - B(\alpha_N)$ since
$Q_N(\alpha_N) =0.$ If $B(\alpha_N) \neq 0$, which is a generic
case, this determines $T(\lambda)$ uniquely (and thus the matrix $M$) by a
triple $(A, B, C).$  This completes the proof of Theorem
1.

{\bf Example.} Consider the simplest case $N=1,\, g=0.$ Then the
$\alpha$-triples have the form $$A = 1, \quad B = \beta, \quad C =
- (\lambda + \gamma)$$ with arbitrary $\beta,\, \gamma \in
\mathbb{C},$ so that in the periodic case the corresponding $\phi$
have a general form
\begin{equation}
\label{ex1}
 \phi = -\beta +
\sqrt{\lambda + \gamma}.
\end{equation}
 In this particular case this can be
easily seen directly. Indeed $\phi =
[b_0,\overline{b_1}]_{\alpha}$ satisfies the quadratic equation
$$\phi^2 + (b_1^* - b_0)\phi - (\lambda - \alpha_1 + b_0b_1^*) =
0,$$ where $b_1^* = b_1 - b_0.$ Thus to find a periodic continued
$\alpha$-fraction expansion of (\ref{ex1}) one should solve the
system of equations
\begin{equation}
\label{sys} b_1^* - b_0 = 2 \beta, \quad b_0b_1^* = \alpha_1 +
\gamma,
\end{equation}
 which has two solutions $$b_0 = -\beta \pm \sqrt{\beta^2 +
\alpha_1 + \gamma}, \quad b_1^* = \beta \pm \sqrt{\beta^2 +
\alpha_1 + \gamma}.$$ One can easily check that these two
solutions correspond to two solutions of the factorisation problem
from the previous section.

In the pure periodic case by the additional condition (\ref{add})
we have $C = -(\lambda - \alpha_1),$ so $\gamma = -\alpha_1$ and
the general form of $\phi$ is
\begin{equation}
\label{ex2}
 \phi = -\beta +
\sqrt{\lambda - \alpha_1}.
\end{equation} In that case we have also $b_1^* =
b_1 - b_0 =0,$ so the system (\ref{sys}) reduces to just one
equation $b_0 = -2\beta,$ which determines the pure periodic
$\alpha$-fraction expansion of (\ref{ex2}) uniquely in
agreement with our previous consideration.

\section{Action of $ \mathbb{Z}_2 \times
S_N$}

A surprising corollary of Theorem 1 is the invariance of the set of
$N$-periodic $\alpha$-fractions under the permutations
$\sigma \in S_N$ of the set $\alpha:$ $$\sigma(\alpha)_k =
\alpha_{\sigma(k)}.$$ This is not obvious from the very beginning
and in fact is not true in the pure periodic case.

In this section we explain how to use this symmetry to describe
all $2 N!$ periodic $\alpha$-fractions for a given
algebraic function $\phi.$ In fact, the full symmetry group is the
product $G = \mathbb{Z}_2 \oplus S_N$, where $\mathbb{Z}_2$ is
generated by the involution $\varepsilon$ interchanging two different
$\alpha$-fraction expansions with the same order of the parameters
$\alpha_i$.

The action of this group is described by Theorem 2. We are
going to prove it now.

Recall that $\pi \in S_N$ is the permutation, which reverses the
order $1, 2, \dots, N-1,\, N$ to $N,\, N-1, \dots, 2, 1,$ and the
involution $\sigma_k$ swaps $k$ and $k+1$ leaving the rest fixed.

Let us start with the action of $\sigma_k$ first. Let us introduce
(assuming that $b_k \neq 0$)
\begin{equation}
\label{actionsigma2} \tilde b_{k-1} = b_{k-1}+
\frac{\alpha_{k+1}-\alpha_k}{b_k} ,\quad \tilde b_{k+1}= b_{k+1}
-\frac{\alpha_{k+1}-\alpha_k}{b_k}, \quad k=1,\dots, N-1.
\end{equation}
One can check directly the following matrix identity: $$\left [
\begin {array}{cc} b_{k-1} & \lambda- \alpha_k \\1&0 \end{array}
\right ]  \left [ \begin {array}{cc} b_{k} & \lambda- \alpha_{k+1}
\\1&0 \end{array} \right ]  \left [ \begin {array}{cc} b_{k+1} &
\lambda- \alpha_{k+2} \\1&0 \end{array} \right ] = $$ $$\left [
\begin {array}{cc} \tilde b_{k-1} & \lambda- \alpha_{k+1} \\1&0 \end{array}
\right ]  \left [ \begin {array}{cc} b_{k} & \lambda- \alpha_{k}
\\1&0 \end{array} \right ]  \left [ \begin {array}{cc} \tilde b_{k+1} &
\lambda- \alpha_{k+2} \\1&0 \end{array} \right ].$$ Similarly for
$k=N-1$ we have $\left [ \begin {array}{cc} b_{N-2} & \lambda-
\alpha_{N-1} \\1&0 \end{array} \right ]  \left [ \begin
{array}{cc} b_{N-1} & \lambda- \alpha_{N} \\1&0 \end{array} \right
]  \left [ \begin {array}{cc} 1 &b_N-b_0\\ 0&1 \end{array} \right
] =\left [ \begin {array}{cc} \tilde b_{N-2} & \lambda- \alpha_{N}
\\1&0
\end{array} \right ]  \left [ \begin {array}{cc} b_{N-1} &
\lambda- \alpha_{N-1} \\1&0 \end{array} \right ]  \left [ \begin
{array}{cc} 1 &\tilde b_N-b_0 \\ 0&1 \end{array} \right ].$ Taking
into account the results of the previous section we see that the
action of $\sigma_k$ is indeed given by the formula (\ref
{actionsigma}).

To prove the remaining part of Theorem 2 recall that $\phi=[b_0;
\overline {b_1, \dots, b_{N-1},b_N}]_{\alpha}$ is the fixed point
of the fractional linear transformation (\ref{fraclin}), and
therefore it is the fixed point of its inverse, which as one can
easily check is given by $$s^{-1}(\phi)= -b_N + b_0+
\frac{a_N}{-b_{N-1} +\frac{a_{N-1}}{-b_{N-2}+_{\ddots
-b_{1}+\frac{a_1}{-b_0+\phi}}}}.$$ Thus
\begin{equation}
\label{second} \phi=[b_0-b_N; \overline
{-b_{N-1},\dots,-b_1,-b_N}]_{\alpha_N,\dots, \alpha_1,\alpha_0}
\end{equation}
 is
a periodic $\alpha$-fraction corresponding to the
sequence $\pi(\alpha)=\alpha_N,\dots, \alpha_1,\alpha_0.$ Now note
that the trace of the corresponding matrix (\ref{mat2}) is equal to
$P_{N-1}+ Q_{N} = (b_1 + b_2 +\dots +b_N) \lambda^g + \dots$ (in
the notations of the previous section). If we replace $b_1,\dots,
b_n$ by  $-b_{N-1},\dots,-b_1,-b_N$ its highest coefficient
clearly changes sign. This means that the new periodic
$\alpha$-fraction (\ref{second}) corresponds to the
action of the element $\varepsilon \pi \in G.$ Since $\varepsilon \pi$ and $\sigma_k$
generate the group $G$ we have described the full action.

In the pure periodic case because of the additional condition
$C(\alpha_N)=0$ the symmetry group is reduced to $S_{N-1}$,
permuting $\alpha_i$ with $i=1,\dots,N-1.$ This group is generated
by $\sigma_k$ with $k=1, \dots, N-2$ with the action given by the
same formula (\ref{actionsigma}). Theorem 2 is proved.

{\bf Example.} Let $N=3,\, g = 1, \, \alpha = (1,3,4)$ and
$$\phi=\frac{3x-7 + \sqrt{4x^3- 31x^2 + 62x +1}}{2(x-6)}.$$ We
have $A(x)= x-6,  \quad B(x) =-\frac{1}{2}(3x-7),\quad C(x) = -
x^2 +4x -2$ and $R(x) = \frac{1}{4}(4x^3- 31x^2 + 62x
+1)=(\frac{x-7}{2})^2+ (x-1)(x-3)(x-4)$ is $\alpha$-admissible.

 To each permutation of the sequence (1,3,4) we have the following two periodic continued
 $\alpha$-fraction representations of $\phi$:

$$\phi = [1;\overline{-3,1,3}]_{1,3,4}=
[-\frac{1}{5};\overline{-\frac{15}{6},\frac{6}{5},\frac{3}{10}}]_{1,3,4}$$

$$\phi = [1;\overline{-2,1,2}]_{1,4,3}=
[-\frac{1}{5};\overline{-\frac{5}{3},\frac{6}{5},-\frac{8}{15}}]_{1,4,3}$$

$$\phi =
[\frac{1}{3};\overline{-3,\frac{5}{3},\frac{7}{3}}]_{3,1,4}=
[-1;\overline{-\frac{15}{6},2,-\frac{1}{2}}]_{3,1,4}$$

$$\phi =
[\frac{1}{3}; \overline{-\frac{6}{5},\frac{5}{3},\frac{8}{15}}]_{3,4,1}=
[-1;\overline{-1,2,-2}]_{3,4,1}$$

$$\phi =
[-\frac{1}{2};\overline{-\frac{6}{5},\frac{15}{6},-\frac{3}{10}}]_{4,3,1}
=[-2;\overline{-1,3,-3}]_{4,3,1}$$

$$\phi =
[-\frac{1}{2};\overline{-2,\frac{15}{6},\frac{1}{2}}]_{4,1,3}=
[-2;\overline{-\frac{5}{3},3,-\frac{7}{3}}]_{4,1,3}.$$

Once one of them is known we can find the rest using the action of
the group $G$ described above.

\section{Relation with the affine hyperelliptic Jacobi varieties}

In this section we will discuss the geometric aspects of the
periodic $\alpha$-fractions.

Let us first note that strictly speaking a function
$$\phi(\lambda)=\frac{-B(\lambda)+\sqrt{R(\lambda)}}{A(\lambda)}$$
is not a function of $\lambda \in \mathbb{C}$, but a function on
the corresponding hyperelliptic curve $\Gamma(R)$ given by the
equation
\begin{equation}
\label{gammar} \mu^2 = R(\lambda).
\end{equation}
This curve has a natural involution $\tau: (\lambda, \mu)
\rightarrow (\lambda, -\mu),$ interchanging the branches of the
square root.

Let us fix now an $\alpha$-admissible polynomial $R(\lambda)$ and
assume that all its roots are distinct, so that the corresponding
hyperelliptic curve $\Gamma_R$ is non-singular and has genus $g.$
Let us ask the following natural question: which functions on $\Gamma_R$
admit periodic $\alpha$-fractions ? According to Theorem
1 the variety of such functions is birationally equivalent to the
variety $M_R^{\alpha}$ of the triples of polynomials $(A, B, C)$
such that
\begin{equation}
\label{mr} A(\lambda) = \lambda^g + \dots,\,\, C(\lambda) =
-\lambda^{g+1} + \dots,\, \, \deg B \leq g,
\end{equation}
satisfying the relation
\begin{equation}
\label{relat} B^2(\lambda) - A(\lambda)C(\lambda) = R(\lambda)
\end{equation}
with given $R$. The triples satisfying conditions (\ref{mr}) form
an affine space $\mathbb{C}^{3g+2},$ while the relation
(\ref{relat}) is equivalent to $2g+1$ algebraic equations on the
coefficients of $A,\, B,\, C.$ Thus $M_R^{\alpha}$ is an affine
algebraic variety.

We claim that $M_R^{\alpha}$ is nothing else but the product
$J(\Gamma_R)^{aff} \times \mathbb{C},$ where $J(\Gamma_R)^{aff}$
is the affine part of the Jacobi variety of the corresponding
hyperelliptic curve (\ref{gammar}). The proof is essentially due
to Jacobi \cite{Jacobi}, who found the following elementary
description of the hyperelliptic Jacobi variety. We follow here
Mumford's lectures \cite{Mum}.

The affine Jacobi variety $J(\Gamma_R)^{aff}$ is defined as the
set of the positive divisors $D = P_1 + \dots + P_g,$ where $P_1,
\dots, P_g$ are points on the affine curve (\ref{gammar}) (not
necessarily distinct) such that $P_i \neq \tau(P_j)$ for any $i
\neq j.$

Consider the following {\it Jacobi triples} $(U,V,W)$ of
polynomials in $\lambda$, where $U$ and $W$ are monic polynomials
of degree $g$ and $g+1$ respectively, $V$ is a polynomial of
degree less or equal than $g-1$, such that the following relation
is satisfied
\begin{equation}
\label{relatj} V^2(\lambda) + U(\lambda)W(\lambda) = R(\lambda).
\end{equation}
We denote the corresponding variety $N_R.$

{\bf Theorem (Jacobi).} {\it There is a natural bijection between
the set $N_R$ of the Jacobi triples and the affine Jacobi variety
$J(\Gamma_R)^{aff}$.}

We sketch here the proof, which is actually not difficult. Let $D
= P_1 + \dots + P_g$ be a divisor from $J(\Gamma_R)^{aff}.$ Let
us assume at the beginning that all points $P_i \in \Gamma_R$ are
distinct. We would like to associate with $D$ a Jacobi triple
$(U,V,W).$  This can be done as follows.

Let $(\mu_i, \lambda_i)$ be coordinates of the points $P_i.$ Then
the polynomial $U$ is defined as
$$U(\lambda)=\prod_{i=1}^{g}(\lambda-\lambda_i).$$ The polynomial
$V$ is defined by the Lagrange interpolation formula from the
conditions $V(\lambda_i) = \mu_i,\, i = 1, \dots, g$ as
$$V(\lambda)=\sum_{i=1}^{g}\mu_i \frac{\prod_{j\neq i}
(\lambda-\lambda_j)}{\prod_{j\neq i} (\lambda_i-\lambda_j)}$$
(note that by assumptions all $\lambda_i$ are distinct). Now $W$
is defined uniquely by the relation (\ref{relatj}). It is a
polynomial because by  construction $R - V^2$ vanishes at all
the zeroes of $U(\lambda).$ If some of the points $P_i$ collide
one should use a natural generalisation involving also the
derivatives (see \cite{Mum} for details). This defines a natural
map from $J(\Gamma_R)^{aff}$ to $N_R$.

The inverse map is also quite natural: the $\lambda$-coordinates
of $P_i$ are given by the zeroes of $U(\lambda)$ while the
$\mu$-coordinates are the values of the polynomial $V$ at these
zeroes. The fact that the corresponding points $P_i$ belong to
$\Gamma_R$ follows from the relation (\ref{relatj}).

Now let us return to our set $M^{\alpha}_R.$ To every Jacobi
triple $(U,V,W)$ and complex number $\beta \in \mathbb{C}$ we can
relate a unique $\alpha$-triple $(A, B, C)$ defined by
\begin{equation}
\label{f} A=U,\quad B = V+\beta U, \quad C = - W + 2\beta V +
\beta^2 U.
\end{equation}
 Indeed $$B^2 - AC = (V+ \beta U)^2 -U(-W + 2\beta V +
\beta^2 U) = V^2 + UW = R.$$

Conversely, for a given $\alpha$-triple $(A, B, C)$ there exists a
unique pair $((U,V,W), \beta)$  where $\beta$ is the coefficient
of $\lambda^g$ in $B$ and

\begin{equation}
\label{finv}U=A,\, V= B -\beta A,\, W = - C + 2\beta B - \beta^2 A
\end{equation}
is a Jacobi triple.

In the pure periodic case we have a well-defined $2:1$ map $f:
P^{\alpha}_R \rightarrow J(\Gamma_R)^{aff}$ defined by the
same formulae (\ref{f}). Indeed for a given Jacobi triple $(U,V,W)$
there are generically 2 $\alpha$-triples $(A, B, C)$ given by
(\ref{finv}), where $\beta$ must satisfy the quadratic equation $$
U(\alpha_N) \beta^2 + 2 V(\alpha_N) \beta - W(\alpha_N) =
C(\alpha_N) =0,$$ so that $\beta=-\frac{V(\alpha_N)\pm
\sqrt{R(\alpha_N)}}{U(\alpha_N)}$ if $U(\alpha_N) \neq 0$ or
$\beta =\frac{W(\alpha_N)}{2V(\alpha_N)}$ if $U(\alpha_N) = 0$ and
$V(\alpha_N) \neq 0.$ Note that if both $U(\alpha_N) = 0$ and
$V(\alpha_N) = 0,$ then $R(\alpha_N)=0$ which we have excluded.

In other words, to any divisor of degree $g$ of the form $D=
\sum_{i=1}^{g} P_i ,$ we construct the polynomial $A(\lambda)$ in
the same way as $U(x).$ The condition $C(\alpha_N)=0$ implies
$B(\alpha_N)=\sqrt{R(\alpha_N)},$ which means that the polynomial
$B(\lambda)$ is now required to pass through the $g$ points $P_i$
plus the additional point $P(\alpha_N, \sqrt{R(\alpha_N)}) \in
\Gamma_R.$ Generically there are two such points depending on the
choice of the square root, which determines the corresponding
$B(\lambda)$ and thus the $\alpha$-triple uniquely.

Conversely, to any  $\alpha$-triple $(A,B,C)$ with $C(\alpha_N)=0$
we can associate a divisor $D$ and a point $P \in \Gamma_R$ such
that $\Lambda(P) = \alpha_N,$ where $\Lambda: \Gamma_R
\rightarrow \mathbb{C}$ be the projection defined by
$\Lambda(\lambda, \mu) = \lambda.$

Let us discuss now the corresponding meromorphic functions $\phi$ on $\Gamma_R.$

Let $D= P_1 + \dots P_g, \, P_i \in \Gamma_R^{aff}$ be a positive non-special divisor, $L(D + P_{\infty})$
be the corresponding linear space of meromorphic functions with poles at most at $P_1, \dots, P_g$ and $P_{\infty}.$  It has dimension 2. Indeed by the classical {\it Riemann-Roch theorem} \cite{GH}
$$
\dim L(D + P_{\infty}) =  2 + \dim \Omega(D+P_{\infty}),
$$ 
where $\Omega(D+P_{\infty})$ is the linear space of holomorphic 1-forms $\omega$ on $\Gamma_R$ with zeroes at $P_1, \dots, P_g$ and $P_{\infty}.$ Any holomorphic 1-form on the hyperelliptic curve $\Gamma_R$ has a form
$$\omega = \frac{P(\lambda) d\lambda}{\sqrt{R(\lambda)}}$$
where $P(\lambda)$ is a polynomial of degree less than or equal to $g-1.$ Since $D$ is non-special,
such $\omega$ must be identically zero so that $\dim \Omega(D+P_{\infty}) = 0$ and $\dim L(D + P_{\infty}) =  2.$

Consider the functions from the linear space  $L(D + P_{\infty})$ which are equivalent to
$\sqrt{\lambda}$ at infinity; these functions differ from each other by additive constant.
One can easily see that up to this freedom they have the form (\ref{Algebraic_fct}) given by the Jacobi construction above. Note the formulae (\ref{f}) corresponds to the shift $\phi \rightarrow \phi + \beta.$
In the pure periodic case  $\phi$ has a zero at one of two points of $\Gamma_R$ with $\lambda= \alpha_N,$ which reduces the shift to two values. This completes the proof of Theorem 3.

\section*{Acknowledgements} 

We are very grateful to Vassilis Papageorgiou
for the numerous stimulating discussions.

This work has been partially supported by the EPSRC and by the
European Union through the FP6 Marie Curie RTN ENIGMA (Contract
number MRTN-CT-2004-5652) and the ESF programme MISGAM.

\end{document}